\documentclass[12pt,reqno]{amsart}

\usepackage{color}
\usepackage{amsmath, amsthm, amssymb}
\usepackage{amsfonts}
\usepackage[ansinew]{inputenc}
\usepackage[dvips]{epsfig}
\usepackage{graphicx}
\usepackage[english]{babel}
\usepackage{hyperref}
\theoremstyle{plain}
\newtheorem{thm}{Theorem}[section]

\theoremstyle{definition}

\theoremstyle{remark}
\newtheorem{rem}[thm]{Remark}

\numberwithin{equation}{section}

\newcommand{\average}{{\mathchoice {\kern1ex\vcenter{\hrule height.4pt
width 6pt depth0pt} \kern-9.7pt} {\kern1ex\vcenter{\hrule
height.4pt width 4.3pt depth0pt} \kern-7pt} {} {} }}
\newcommand{\ave}{\average\int}
\def\R{\mathbb{R}}

\begin{document}

\title[Obstacle problems and free boundaries: an overview]{Obstacle problems and free boundaries: \\ an overview}

\author{Xavier Ros-Oton}
\address{The University of Texas at Austin, Department of Mathematics, 2515 Speedway, Austin, TX 78751, USA}
\email{ros.oton@math.utexas.edu}


\thanks{This paper originates from my lecture at the XXV CEDYA conference, where I was awarded the 2017 Antonio Valle Prize from the Sociedad Espa\~nola de Matem\'atica Aplicada (SeMA). \\I gratefully thank the SeMA and the organizers of the conference for such opportunity.}


\maketitle

\begin{abstract}
Free boundary problems are those described by PDEs that exhibit a priori unknown (free) interfaces or boundaries.
These problems appear in Physics, Probability, Biology, Finance, or Industry, and the study of solutions and free boundaries uses methods from PDEs, Calculus of Variations, Geometric Measure Theory, and Harmonic Analysis.
The most important mathematical challenge in this context is to understand the structure and regularity of free boundaries.

In this paper we provide an invitation to this area of research by presenting, in a completely non-technical manner, some classical results as well as some recent results of the author.
\end{abstract}

\vspace{5mm}

\section{Introduction}

Many problems in Physics, Industry, Finance, Biology, and other areas can be described by PDEs that exhibit apriori unknown (free) interfaces or boundaries.
These are called Free Boundary Problems.

A classical example is the Stefan problem, which dates back to the 19th century \cite{Stefan,Stefan2}.
It describes the melting of a block of ice submerged in liquid water.
In the simplest case (the one-phase problem), there is a region where the temperature is positive (liquid water) and a region where the temperature is zero (the ice).
In the former region the temperature function $\theta(t,x)$ solves the heat equation (i.e., $\theta_t=\Delta \theta$ in $\{(t,x)\,:\,\theta(t,x)>0\}$), while in the other region the temperature $\theta$ is just zero.
The position of the free boundary that separates the two regions is part of the problem, and is determined by an extra boundary condition on such interface (namely, $|\nabla_x \theta|^2=\theta_t$ on $\partial\{\theta>0\}$).

Another classical and important free boundary problem is the {obstacle problem}.
The solution $u(x)$ of the problem can be thought as the equilibrium position of an elastic membrane whose boundary is held fixed, and which is constrained to lie above a given obstacle $\varphi(x)$ (see Figure 1).
\begin{figure}[htp]
\begin{center}
\includegraphics[]{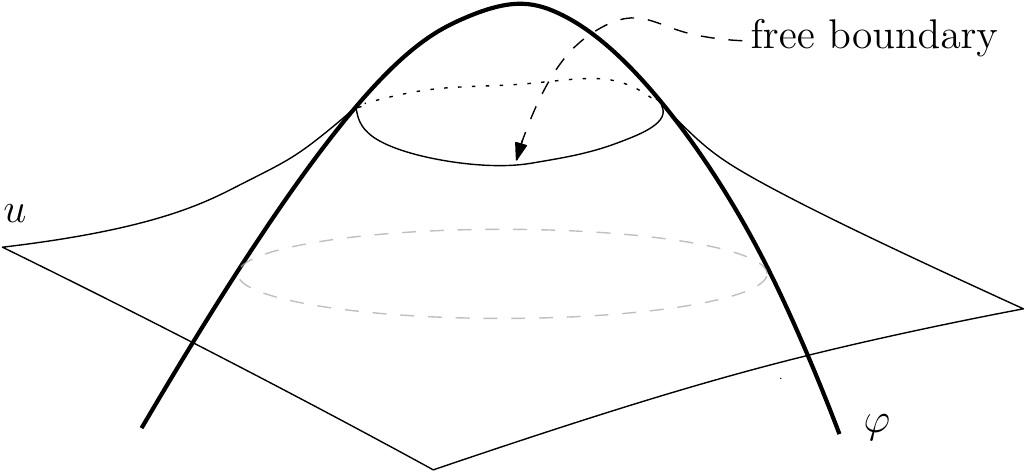}
\end{center}
\caption{The obstacle $\varphi$, the solution $u$, and the free boundary $\partial\{u>\varphi\}$.}
\end{figure}
In the region where the membrane is above the obstacle, the solution $u$ solves an elliptic PDE (say, $\Delta u=0$ in $\{u>\varphi\}$), while in the other region the membrane coincides with the obstacle (i.e., $u=\varphi$).
As in the Stefan problem, there is an extra condition that makes the problem well-posed and determines indirectly the position of the free boundary (in this case, $-\Delta u\geq0$ across the free boundary $\partial\{u>\varphi\}$).

It turns out that, after an appropriate transformation (see \cite{Duv}), the Stefan problem described above is exactly equivalent to the parabolic obstacle problem (in which the ``membrane'' $u$ would evolve with time).
Thus, both problems lie in the setting of \emph{obstacle problems}.
This type of problems appear in fluid mechanics, elasticity, probability, finance, biology, and industry.

From the mathematical point of view, the most challenging question in these problems is to understand the \emph{regularity of free boundaries}.
For example, in the Stefan problem: is there a regularization mechanism that smoothes out the free boundary, independently of the initial data? (Notice that a priori the free boundary could be a very irregular set, even a fractal set!)
Such type of questions are usually very hard, and even in the simplest cases almost nothing was known before the 1970s.
The development of the regularity theory for free boundaries started in the late seventies, and since then it has been a very active area of research.

Generally speaking, a key difficulty in the study of the regularity of free boundaries is the following.
For solutions to elliptic (or parabolic) PDE, one has an equation for a function $u$, and such equation forces $u$ to be regular.
For example, for harmonic functions (i.e., $\Delta u=0$) the equation yields the mean value property, which in turn implies that $u$ is smooth.
In free boundary problems such task is much more difficult.
One does not have a regularizing equation for the free boundary, but only an equation for $u$ which indirectly determines the free boundary.
The classical regularity theory for elliptic/parabolic PDE does not apply, and in many cases one needs to combine techniques from PDEs with ideas and tools from Geometric Measure Theory, Calculus of Variations, and Harmonic Analysis; see \cite{PSU}.

\vspace{3mm}

In this paper, we briefly describe the classical regularity theory for free boundaries in obstacle problems \cite{C-obst1,ACS,CSS}, as well as some new results recently obtained by the author \cite{CRS,BFR2}.
We will study the following three problems:
\begin{itemize}
\item The classical obstacle problem
\item The thin obstacle problem
\item Obstacle problems for integro-differential operators
\end{itemize}
We first present in Section \ref{sec2} the mathematical formulation of these three problems.
Then, in Section \ref{sec3} we explain some motivations and applications.
Finally, in Sections \ref{sec4} and \ref{sec5} we describe the classical regularity theory for free boundaries in these problems, as well as our recent results.

\section{Obstacle problems}
\label{sec2}

\subsection{The classical obstacle problem}

The obstacle problem is probably the most classic and motivating example in the study of variational inequalities and free boundary problems.
Its simplest mathematical formulation is to seek for minimizers of the Dirichlet energy functional
\begin{equation}\label{Dir}
\mathcal E(u)=\int_D |\nabla u|^2dx
\end{equation}
among all functions $u$ satisfying $u\geq\varphi$ in $D$, for a given smooth obstacle $\varphi\in C^\infty$.
Here, $D\subset \R^n$, and one usually has Dirichlet boundary conditions $u=g$ on $\partial D$.
When $D=\R^n$, one simply prescribes $u\to0$ at $\infty$.

A simple variational argument shows that the Euler-Lagrange equations of such minimization problem are
\begin{eqnarray}
\nonumber u&\geq& \varphi \quad \textrm{ in}\ D\\
\label{obst} \Delta u&=&0 \quad\, \textrm{ in}\ \{u>\varphi\}\\
\nonumber -\Delta u&\geq& 0\quad\, \textrm{ in}\ D.
\end{eqnarray}
In other words, the solution $u$ is above the obstacle $\varphi$, it is harmonic whenever it does not touch the obstacle, and moreover it is superharmonic everywhere.

The domain $D$ will be split into two regions: one in which the solution $u$ is harmonic, and one in which the solution equals the obstacle.
The latter region is known as the \emph{contact set} $\{u=\varphi\}$.
The interface that separates these two regions is the \emph{free boundary}.

\begin{figure}[htp]
\begin{center}
\includegraphics[]{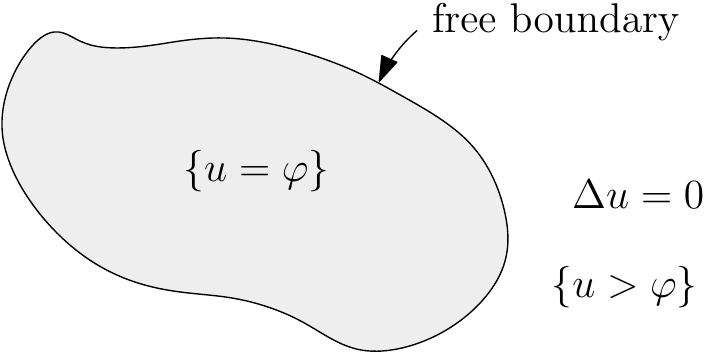}
\end{center}
\caption{The contact set and the free boundary in the classical obstacle problem.}
\end{figure}

Existence and uniqueness of solutions to such problem is not difficult to show.
Indeed, notice that
\[K:=\bigl\{u\in H^1(D)\,:\,u|_{\partial D}=g\quad \textrm{and}\quad u\geq\varphi\bigr\}\]
is a closed and convex subset of $H^1(D)$, and thus there is a unique function $u\in K$ that minimizes the Dirichlet energy \eqref{Dir}.

An alternative way to construct solutions is by using the theory of viscosity solutions.
Namely, the solution of the obstacle problem can be constructed as the least supersolution that lies above the obstacle $\varphi$.
With such approach, a natural way to write the equations \eqref{obst} is
\begin{equation}\label{obst-min}
\min\bigl\{-\Delta u,\,u-\varphi\bigr\}=0\quad \textrm{ in}\quad D.
\end{equation}
Finally, an important remark is that, after the transformation $u\mapsto u-\varphi$, the obstacle problem  \eqref{obst} is equivalent to
\begin{eqnarray}
\nonumber u&\geq& 0 \qquad\ \, \textrm{ in}\ D\\
\label{obstB} -\Delta u&=&f(x) \quad\, \textrm{ in}\ \{u>0\}\\
\nonumber -\Delta u&\geq& f(x)\quad\, \textrm{ in}\ D,
\end{eqnarray}
where $f:=\Delta\varphi$.

It is not difficult to show that solutions to the obstacle problem are $C^{1,1}$ (second derivatives of $u$ are bounded but not continuous).
The most important mathematical challenge in the study of such problem is to understand the regularity of the free boundary.
At first glance, it is not clear at all why should we expect the free boundary to be smooth, since a priori the equations \eqref{obst} do not seem to force it to be regular.
As we will see later, it turns out that there is a regularization effect on the free boundary, and that it is $C^\infty$ (maybe outside a certain set of singular points), as long as the obstacle $\varphi$ is $C^\infty$.

\subsection{The thin obstacle problem}

The thin obstacle problem (also called the boundary obstacle problem) arises when minimizing the Dirichlet energy
\[\mathcal E(u)=\int_{D^+}|\nabla u|^2dx\]
among all functions $u$ satisfying
\[u\geq \varphi\qquad \emph{on} \quad\{x_n=0\}\cap D^+.\]
Here, $D^+\subset B_1\cap \{x_n\geq0\}$, and usually one would take either $D^+=B_1^+$ or $D^+=\R^n_+$.
When $D^+$ is bounded, then the Dirichlet boundary conditions are $u=g$ on $\partial D^+\cap \{x_n>0\}$, while when $D^+=\R^n_+$ one simply prescribes $u\to0$ at $\infty$.

A simple variational argument shows that the Euler-Lagrange equations of such minimization problem are
\begin{equation}\begin{split}\label{thinobst}
 \Delta u&=0 \quad\, \textrm{ in}\ D^+\cap \{x_n>0\}\\
  u&\geq \varphi \quad \textrm{ in}\ D^+\cap \{x_n=0\}\\
 \partial_{x_n}u&\leq 0\quad\, \textrm{ in}\ D^+\cap \{x_n=0\}\\
 \partial_{x_n}u&= 0\quad\, \textrm{ in}\ D^+\cap \{x_n=0\}\cap \{u>\varphi\}.
\end{split}\end{equation}
As in the classical obstacle problem, the existence and uniqueness of solutions for such problem is standard.

\begin{figure}[htp]
\begin{center}
\includegraphics[]{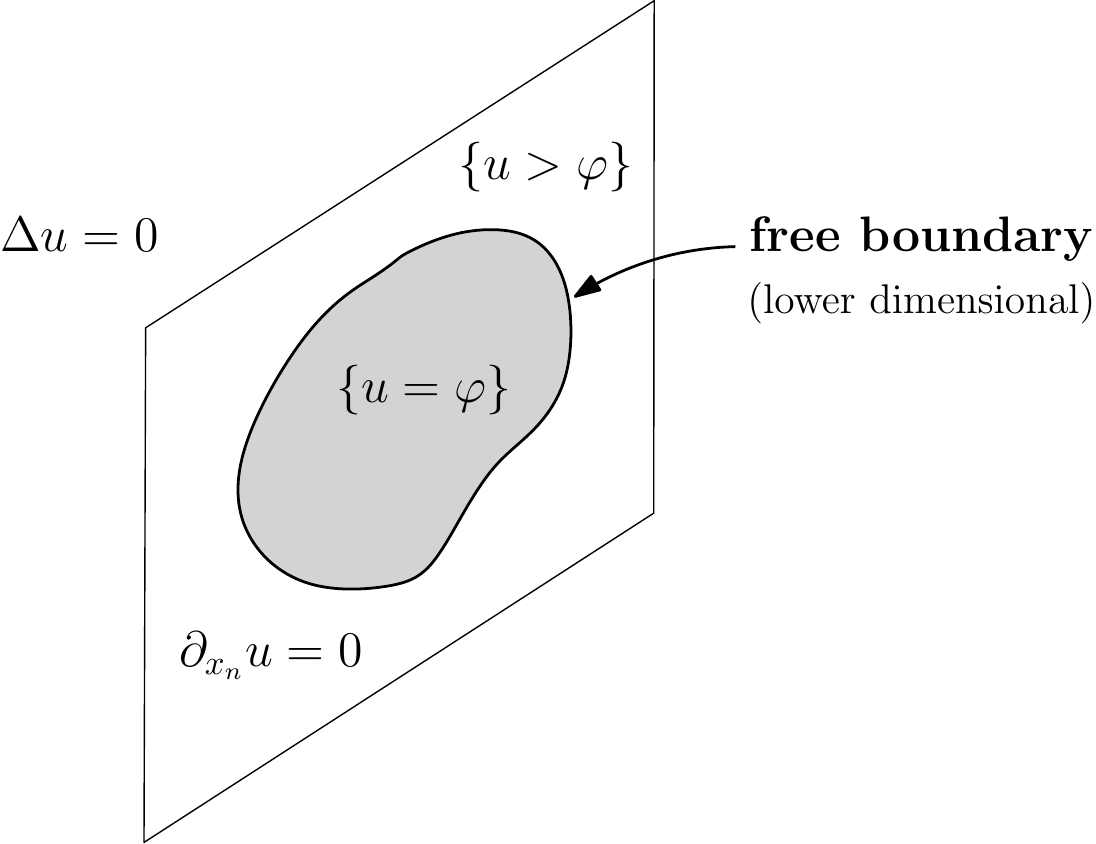}
\end{center}
\caption{The contact set and the free boundary in the thin obstacle problem.}
\end{figure}

The set $D^+\cap \{x_n=0\}$ will be split into two regions: one in which $\partial_{x_n}u$ is zero, and one in which $u$ equals the obstacle.
The latter region is the \emph{contact set}.
The interface that separates these two regions is the \emph{free boundary}.

The equations \eqref{thinobst} can be written as
\begin{equation}\label{thinobst-min}
\begin{split}
\Delta u&=0\quad \textrm{ in}\ D^+\cap \{x_n>0\}\\
\min\bigl\{-\partial_{x_n} u,\,u-\varphi\bigr\}&=0\quad \textrm{ on}\ D^+\cap \{x_n=0\},
\end{split} \end{equation}
which is the analogous of \eqref{obst-min}.

After an even reflection with respect to the hyperplane $\{x_n=0\}$, the solution $u$ will be harmonic across such hyperplane wherever $\partial_{x_n}u=0$, and it will be superharmonic wherever $\partial_{x_n}u<0$.
Thus, such reflected function would formally solve the classical obstacle problem, but with the restriction $u\geq\varphi$ only on $\{x_n=0\}$ (the obstacle is thin).

\subsection{Obstacle problems for integro-differential operators}

A more general class of obstacle problems is obtained when minimizing \emph{nonlocal} energy functionals of the form
\[\mathcal E(u)=\int_{\R^n}\int_{\R^n}\bigl|u(x)-u(y)\bigr|^2K(x-y)\,dx\,dy\]
among all functions $u\geq\varphi$ in $\R^n$ ---or with $u=g$ in $D^c$ and $u\geq\varphi$ in $D$.
Here, $K$ is a nonnegative and even kernel ($K\geq0$ and $K(z)=K(-z)$), and the minimal integrability assumption is
\[\int_{\R^n}\min\bigl\{1,\,|z|^2\bigr\}K(z)dz<\infty.\]
The most simple and canonical example is
\begin{equation}\label{power}
\qquad \qquad K(z)=\frac{c}{|z|^{n+2s}},\qquad \qquad s\in(0,1),
\end{equation}
while a typical ``uniform ellipticity'' assumption is
\begin{equation}\label{L0}
\frac{\lambda}{|z|^{n+2s}}\leq K(z)\leq \frac{\Lambda}{|z|^{n+2s}},
\end{equation}
with $s\in(0,1)$ and $0<\lambda\leq \Lambda$; see for example \cite{BL,CS,Ros-survey}.

The Euler-Lagrange equations of such minimization problem are
\begin{equation}\begin{split}\label{obst-nonl}
u&\,\geq\, \varphi \quad \textrm{ in}\quad \R^n\\
L u&\,=\,0 \quad\, \textrm{ in}\quad \{u>\varphi\}\\
-L u&\,\geq\, 0\quad\, \textrm{ in}\quad \R^n,
\end{split}\end{equation}
where $L$ is an integro-differential operator of the form
\begin{equation}\label{L}\begin{split}
Lu(x)&={\rm p.v.}\int_{\R^n}\bigl(u(y)-u(x)\bigr)K(x-y)dy\\
&={\rm p.v.}\int_{\R^n}\bigl(u(x+z)-u(x)\bigr)K(z)dz\\
&=\frac12\int_{\R^n}\bigl(u(x+z)+u(x-z)-2u(x)\bigr)K(z)dz.
\end{split}\end{equation}
In other words, $u$ solves the obstacle problem \eqref{obst} but with the Laplacian $\Delta$ replaced by the integro-differential operator $L$ in \eqref{L}.

\vspace{3mm}

\begin{rem}[Relation to the thin obstacle problem]\label{rem-thin}
When $K$ is given by \eqref{power}, then $L$ is a multiple of the \emph{fractional Laplacian} $-(-\Delta)^s$, defined by
\[(-\Delta)^s u(x)=c_{n,s}{\rm p.v.}\int_{\R^n}\frac{u(x)-u(y)}{|x-y|^{n+2s}}\,dy.\]
(The constant $c_{n,s}$ is chosen so that the Fourier symbol of this operator is $|\xi|^{2s}$.)
The obstacle problem for the fractional Laplacian extends at the same time the classical obstacle problem and the thin obstacle problem.

Indeed, on the one hand, as $s\uparrow 1$, such operator converges to the Laplacian operator $\Delta$.
On the other hand, when $s=\frac12$ the half-Laplacian $(-\Delta)^{1/2}$ can be written as a Dirichlet-to-Neumann operator in $\R^{n+1}_+$:
for any function $w(x)$ in $\R^n$, if we consider its harmonic extension $\tilde w(x,x_{n+1})$ in $\R^{n+1}_+$ then the Neumann derivative $\partial_{x_{n+1}}\tilde w$ on $\{x_{n+1}=0\}$ is exactly the half-Laplacian of $w(x)$ as a function on $\R^n$.
\begin{figure}[htp]
\begin{center}
\includegraphics[]{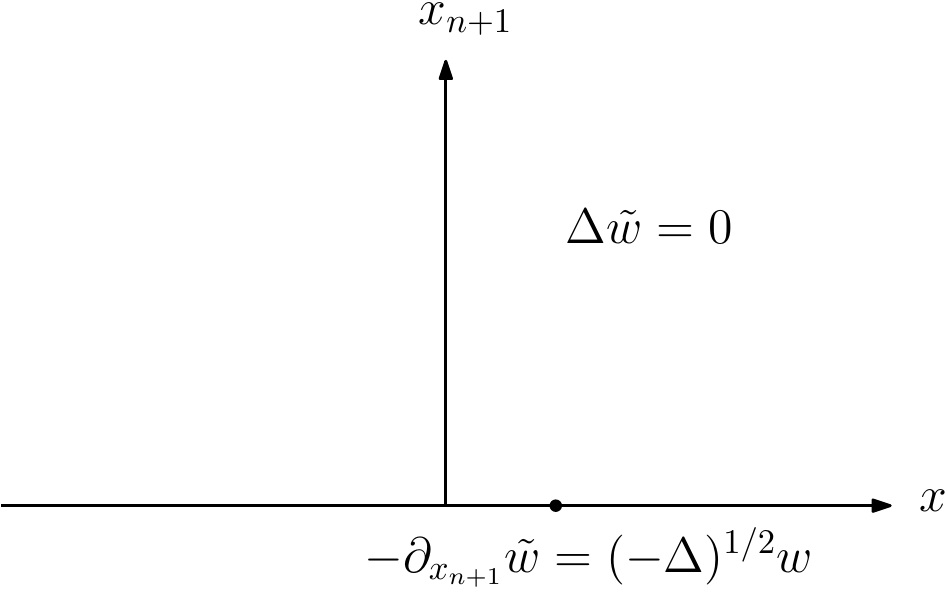}
\end{center}
\caption{The half-Laplacian in $\R^n$ as a Dirichlet-to-Neumann operator in $\R^{n+1}_+$.}
\end{figure}

Therefore, the thin obstacle problem \eqref{thinobst-min} (with $D^+=\R^n_+$) is the same as
\[\min\bigl\{(-\Delta)^{1/2}u,\,u-\varphi\bigr\}=0\quad \textrm{ in}\quad \R^{n-1}.\]
Notice that here we just consider the function $u$ on $\{x_n=0\}$, and this is why the problem is in one dimension less, $\R^{n-1}$.
Note also that with this alternative formulation of the thin obstacle problem the free boundary is not lower-dimensional anymore, but the operator has changed and it is now $(-\Delta)^{1/2}$.

To sum up, the obstacle problem for the fractional Laplacian includes the thin obstacle problem as a particular case ($s=\frac12$), and also includes, as a limiting case $s\uparrow1$, the classical obstacle problem.
\end{rem}

\section{Motivation and applications}
\label{sec3}

A remarkable fact about obstacle problems is that the same problem may appear in many different contexts.
We will next briefly describe some of these applications, and we refer to the books \cite{obst-appl4}, \cite{DL}, \cite{Chipot}, \cite{obst-appl1}, \cite{CT}, \cite{obst-appl2}, and \cite{PSU} for more information.

\subsection{The Stefan problem}

As explained in the Introduction, the Stefan problem describes the temperature distribution $\theta(t,x)$ in a homogeneous medium undergoing a phase change, for example melting ice submerged in liquid water.
In the simplest case, the function $\theta(t,x)$ solves the heat equation $\theta_t-\Delta\theta=0$ in the set $\{\theta>0\}$, and it is just zero otherwise.
The evolution of the free boundary that separates the two phases is determined by the extra boundary condition $|\nabla_x \theta|^2=\theta_t$ on $\partial\{\theta>0\}$.

The problem was formulated by Stefan in \cite{Stefan} and \cite{Stefan3}, but it was first studied many years before by Lam\'e and Clapeyron in \cite{Stefan2}.

It was proved by Duvaut \cite{Duv} that, if one considers the function
\[u(t,x)=\int_0^t \theta(\tau,x)d\tau,\]
then this Stefan problem is transformed into the parabolic obstacle problem
\begin{eqnarray*}
u&\geq& 0 \quad\, \textrm{ in}\ D\\
u_t - \Delta u&=&-1 \quad\, \textrm{ in}\ \{u>0\}\cap D\\
u_t - \Delta u&\geq &-1\quad\, \textrm{ in}\ D;\\
\end{eqnarray*}
which is the evolutionary analogue of \eqref{obstB}.

\subsection{Fluid filtration, constrained heating}

The classical obstacle problem and the thin obstacle problem appear naturally when describing various processes in Mechanics, Biology, or Industry.
Many of these models are described in the books of Rodrigues \cite{obst-appl1}, Duvaut and Lions \cite{DL}, and Chipot \cite{Chipot}, we refer to these books for more details.

\vspace{3mm}

\noindent $\bullet$ A first important example is the \emph{Dam problem}, which describes fluid filtration through porous medium.
The physical problem is the following.
Two water reservoirs, of different levels, are separated by an earth dam: water flows from the highest level to the lowest one, and one looks for the quantities associated to the motion.
It turns out that, after an appropriate transformation (known as the Baiocchi transform  \cite{Baiocchi}), the problem can be reduced to the classical obstacle problem in $\R^3$.

\vspace{3mm}

\noindent $\bullet$ Another important example is the process of \emph{osmosis} in the study of \emph{semipermeable membranes}.
In this setting, one has a membrane and there is fluid flow from one side of the membrane to the other.
The membrane is permeable only for a certain type of molecules (say, water molecules), and blocks other molecules (the solutes).
Water flows from the region of smaller concentration of solute to the region of higher concentration (\emph{osmotic pressure}), and the flow stops when there is enough pressure on the other side of the membrane to compensate for the osmotic pressure.
In the simplest and stationary case, the mathematical formulation of the problem is exactly the thin obstacle problem \eqref{thinobst-min}.

\vspace{3mm}

\noindent $\bullet$ Finally, other examples include models related to interior or boundary \emph{heat control}, which lead to the classical or the thin obstacle problem, respectively.
In such type of models, the temperature function plays essentially an analogous role to that of the pressure in models of Fluid Mechanics; see \cite{DL}.

\subsection{Elasticity}

The classical obstacle problem has a clear interpretation in terms of an elastic membrane which is constrained to be above a given obstacle $\varphi$, and whose boundary is held fixed.
With that interpretation, it is clear that the equilibrium position of the membrane is the minimizer of the energy functional constrained to the condition $u\geq\varphi$.

However, it is maybe the \emph{thin} obstacle problem the one that is more relevant in elasticity.
Indeed, it is very related to the \emph{Signorini problem} \cite{Signorini0,Signorini}, a classical problem in linear elasticity which dates back to 1933.
The problem consists in finding the elastic equilibrium configuration of a 3D elastic body, resting on a rigid frictionless surface and subject only to its mass forces.
The problem leads to a system of variational inequalities for the displacement vector $\vec{u}$ in $\R^3$.
In case that $\vec{u}$ is a scalar function $u$, then the system boils down to the {thin} obstacle problem \eqref{thinobst-min}.

\subsection{Optimal stopping and Mathematical Finance}

Obstacle problems \eqref{obst-min} and \eqref{obst-nonl} arise naturally in Probability and Finance, too.

Indeed, let us consider the following stochastic control model.
We have a stochastic process $X_t$ in $\R^n$ and some given payoff function $\varphi:\R^n\rightarrow\R$.
One wants to find the \emph{optimal stopping} strategy to maximize the expected value of $\varphi$ at the end point.
If we denote $L$ the infinitesimal generator of the process $X_t$, then it turns out that the value function $u(x)$ (i.e., the maximum expected payoff we can obtain starting at $x$) solves the following problem
\begin{equation} \label{obstacle2}
\begin{split}
u\geq\varphi & \quad \textrm{in}\ \R^n,\\
-Lu\geq0 & \quad \textrm{in}\ \R^n,\\
Lu=0 & \quad \textrm{in}\ \{u>\varphi\}.
\end{split}\end{equation}
This means that the value function in any optimal stopping problem solves an obstacle problem.

In the context of financial mathematics this type of problem appears as a model for pricing American options \cite{obst-finance2,obst-appl9}.
The function $\varphi$ represents the option's payoff, and the contact set $\{u=\varphi\}$ is the exercise region.
Notice that, in this context, the most important unknown to understand is the exercise region, i.e., one wants to find and/or understand the two regions $\{u=\varphi\}$ (in which we should exercise the option) and  $\{u>\varphi\}$ (in which we should wait and not exercise the option yet).
The free boundary is the separating interface between these two regions.
In fact, in such models the option usually has an expiration date $T$, and then the obstacle problem to be studied is the \emph{parabolic} version of \eqref{obstacle2}.

\begin{figure}[htp]
\includegraphics{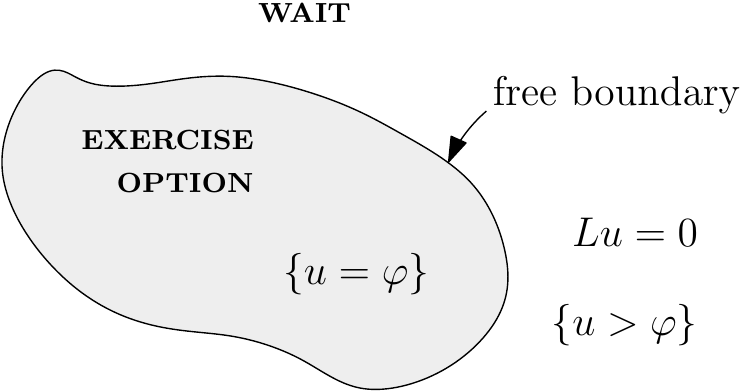}
\caption{The free boundary separates the two regions: the one in which we should exercise the option, and the one in which it is better to wait.}
\end{figure}

If the stochastic process is the Brownian motion, then $L=\Delta$ and the function $u$ will satisfy the classical obstacle problem.
However, in finance one usually needs to consider \emph{jump processes} \cite{obst-appl6,obst-finance,obst-finance2}, and then the operator $L$ will be an integro-differential operator of the form
\[Lu(x)=\,\textrm{p.v.}\int_{\R^n}\bigl(u(y)-u(x)\bigr)K(x-y)\,dy.\]
Such type of operators were introduced in Finance in the seventies by the Nobel Prize winner R. Merton; see \cite{Merton}.

We refer to the book \cite{CT} for a detailed description of the model; see also \cite{obst-finance2,obst-finance,obst-appl6,obst-appl5,obst-appl9}.


\subsection{Interacting energies in physical, biological, or material sciences}

Many different phenomena in Physics, Biology, or material science give rise to models with interacting particles or individuals.

In such context, the mathematical model is usually the following (see, e. g.,~\cite{obst-appl7}).
We are given a repulsive-attractive interaction potential $W\in L^1_{\rm loc}$ in $\R^2$, and the associated interaction energy
\[E[\mu]:=\frac12\int_{\R^2}\int_{\R^2}W(x-y)d\mu(x)d\mu(y),\]
where $\mu$ is any probability measure.

The potential $W$ is repulsive when the particles or individuals are very close, and attractive when they are far from each other.
A typical assumption is that near the origin we have
\begin{equation}\label{W}
W(z)\sim -\frac{1}{|z|^{\beta}}\quad \textrm{for}\quad z\sim 0,
\end{equation}
for some $\beta\in(0,2)$.
It would usually grow at infinity, say $W(z)\sim |z|^{\gamma}$ for $z\sim \infty$.

An important question to be understood is that of the \emph{regularity of minimizers}, i.e., the regularity properties of the measures $\mu_0$ which minimize the interaction energy $E[\mu]$.
For this, it was rigorously shown in \cite{CDM} that the minimizer $\mu_0$ is given by $\mu_0=Lu$, where $u(x):=\int_{\R^n} W(x-y)d\mu_0(y)$ satisfies the \emph{obstacle problem}
\[\min\bigl\{-Lu,\,u-\varphi\bigr\}=0,\]
for a certain operator $L$ and a certain obstacle $\varphi$ that depend on $W$.
When $W$ satisfies \eqref{W} (as well as some extra conditions at infinity), such operator $L$ turns out to be an \emph{integro-differential} operator of the form \eqref{L}, with
\[K(z)\sim \frac{1}{|z|^{4-\beta}}\quad \textrm{for}\quad z\sim 0.\]
In the simplest case, we would have \eqref{power} (with $n=2$ and $2s=2-\beta$).

In $\R^3$, the Newtonian potential $W(z)=|z|^{-1}$ leads to the classical obstacle problem, i.e., we get $L=\Delta$.
When $W$ is more singular than the Newtonian potential ---i.e., $\beta>1$ in \eqref{W}--- then we get an integro-differential operator of the form \eqref{L} satisfying~\eqref{L0}.

Summarizing, understanding the regularity of minimizers $\mu_0$ of the interaction energy $E[\mu]$ is equivalent to understanding the regularity of solutions and free boundaries in obstacle problems for integro-differential operators.
In this setting, the contact set $\{u=\varphi\}$ is the support of the minimizer $\mu_0$.

\section{Regularity theory: known results}
\label{sec4}

Let us next describe the main known mathematical results on the classical obstacle problem, the thin obstacle problem, and obstacle problems for integro-differential operators.

As explained in Section \ref{sec2}, the existence and uniqueness of solutions follows by standard techniques: the solution can be constructed either by minimizing an energy functional among all functions satisfying $u\geq\varphi$, or by using the theory of viscosity solutions (see e.g. \cite{LS,Chipot}).

The central mathematical challenge in obstacle problems is to \emph{understand the geometry and regularity of the free boundary}, i.e., of the interface $\partial\{u>\varphi\}$.
A priori such interface could be a very irregular object, even a fractal set with infinite perimeter.
As we will see, it turns out that this cannot happen, and that free boundaries are smooth (maybe outside a certain set of singular points).

Our presentation will be very brief and will describe only some of the main ideas in the proofs.
We refer the interested reader to the book \cite{PSU} and to the survey \cite{CS-heur} for more details and references.

\subsection{The classical obstacle problem}

The regularity theory for free boundaries in the classical obstacle problem was mainly developed in the seventies, with the groundbreaking paper of L. Caffarelli \cite{C-obst1}.

Namely, the first results for this problem established the optimal $C^{1,1}$ regularity of solutions (i.e., second derivatives of $u$ are bounded but not continuous).
Then, the first general result for free boundaries was proved by Kinderlehrer and Nirenberg \cite{KN}, who showed that, if the free boundary is $C^1$, then it is $C^\infty$.
This is a perturbative result that is proved by flattening the (free) boundary ---via a hodograph transform--- and then using a bootstrap argument.
The main open problem was still open: to understand what happens in general with the regularity of the free boundary.
As said before, a priori it could be a very irregular set with infinite perimeter, while in order to apply the results of \cite{KN} one needs the free boundary to be at least $C^1$.
The breakthrough came with the work \cite{C-obst1}, where Caffarelli developed a regularity theory for free boundaries in the obstacle problem, and established the regularity of free boundaries near regular points.
After that, the set of singular points was studied in dimension $n=2$ in \cite{CR}, and then in higher dimensions in \cite{C-obst2,Mon}.

\begin{figure}[htp]
\begin{center}
\includegraphics{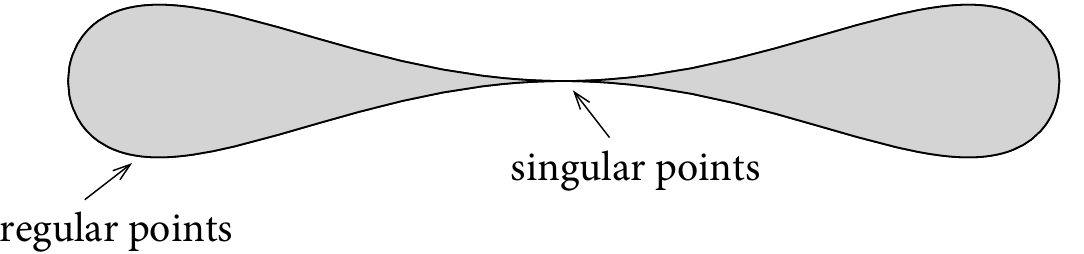}
\end{center}
\caption{A free boundary with a singular point. The contact set $\{u=\varphi\}$ (colored gray) has zero density at the singular point.}
\end{figure}

The main known results from \cite{C-obst1,KN,CR,C-obst2,Mon} can be summarized as follows:
\begin{itemize}
\item[-] At every free boundary point $x_0$ one has
\begin{equation}\label{nondeg}
0<cr^2\leq \sup_{B_r(x_0)}(u-\varphi)\leq Cr^2\qquad r\in(0,1)
\end{equation}

\item[-] The free boundary splits into {regular} points and {singular} points

\item[-] The set of {regular points} is an open subset of the free boundary and it is $C^\infty$

\item[-] {Singular points} are those at which the contact set $\{u=\varphi\}$ has density zero, and these points (if any) are locally contained in a $(n-1)$-dimensional $C^1$ manifold
\end{itemize}

\vspace{2mm}

\noindent Summarizing, \emph{the free boundary is smooth, possibly outside certain set of singular points}.

\vspace{2mm}

To prove such regularity results, one considers \emph{blow-ups}.
Namely, given a free boundary point $x_0$ one shows that
\[u_r(x):=\frac{(u-\varphi)(x_0+rx)}{r^2} \ \longrightarrow\  u_0(x)\qquad \textrm{in}\ \,C^1_{\rm loc}(\R^n),\]
for some function $u_0$ which is a global solution of the obstacle problem.
Notice that the rescaling parameter $r^2$ comes from \eqref{nondeg}.

Then, the main difficulty is to \emph{classify blow-ups}, i.e., show that \vspace{1mm}
\[\textrm{regular point}\quad \Longrightarrow\quad u_0(x)=(x\cdot e)_+^2\quad \textrm{(1D solution)}\]
\[\textrm{singular point}\quad \Longrightarrow\quad u_0(x)=\sum_i \lambda_i x_i^2\quad \textrm{(paraboloid)};\]
see Figure~7.
Notice that, after the blow-up, the contact set $\{u_0=0\}$ becomes a half-space in case of regular points, while it has zero measure in case of singular points.

\begin{figure}[htp]
\begin{center}
\includegraphics{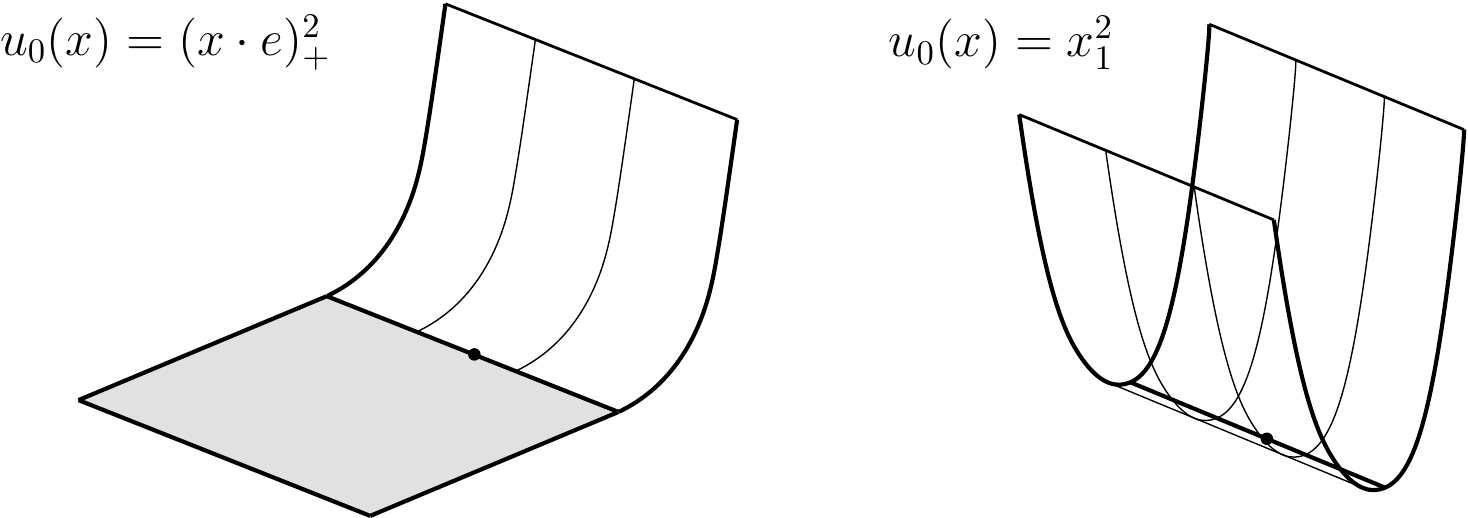}
\end{center}\label{blow-ups}
\caption{The blow-up profile $u_0$ at a regular point (left) and at a singular point (right).}
\end{figure}

Finally, once the classification of blow-ups is well understood, then one has to transfer the information from $u_0$ to $u$, and show that if $x_0$ was a regular point, then the free boundary is $C^1$ in a neighborhood of~$x_0$.
We refer to \cite{C-obst2} or \cite{PSU} for more details.

\subsection{The thin obstacle problem}

The regularity theory for the free boundary differs substantially if we consider the \emph{thin} obstacle problem \eqref{thinobst-min} instead of the classical one \eqref{obst-min}.

In the classical obstacle problem, all blow-ups are homogeneous of degree 2, and the full structure of the free boundary is completely understood, as explained above.
In the thin obstacle problem, instead, understanding the regularity of free boundaries is much harder.
An important difficulty comes from the fact that in thin obstacle problems there is no a priori analogous of \eqref{nondeg}, and thus blow-ups may have different homogeneities.

The first results for the thin obstacle problem were obtained in the 1960's and 1970's.
However, even if the regularity of free boundaries in the classical obstacle problem had been established in 1977 \cite{C-obst1}, nothing was known for the thin obstacle problem.
Such question remained open for 30 years, and was finally answered by Athanasopoulos, Caffarelli, and Salsa in \cite{ACS}.

The main result in \cite{ACS} establishes that if $u$ solves the thin obstacle problem \eqref{thinobst-min} with $\varphi\equiv0$, then for every free boundary point $x_0$ we have
\begin{itemize}
\item[(i)] either \vspace{-4mm}
    \[\qquad \qquad \qquad 0<cr^{3/2}\leq \sup_{B_r(x_0)} u\leq Cr^{3/2}\qquad \textrm{(regular points)}\]
\item[(ii)] or \hspace{3.9cm} $\displaystyle 0\leq \sup_{B_r(x_0)} u\leq Cr^2$.
\end{itemize}
Moreover, they proved that set of {regular points} (i) is an open subset of the free boundary, and it is $C^{1,\alpha}$ for some small $\alpha>0$.

\vspace{2mm}

The proof of this result is strongly related to the theory of \emph{minimal surfaces}; see the survey paper \cite{CS-heur}.
Namely, to study the regularity of the free boundary they found a quantity that is monotone as we zoom in a solution at a given free boundary point.
Generally speaking, {monotonicity formulas} are a kind of ``radial entropy'' that increase as we zoom out at a given point.
This important property usually yields that the blow-up (of a minimal surface, or of a solution to a PDE) has a special configuration.
In the theory of minimal surfaces, the corresponding monotonicity formula implies that the blow-ups of a minimal surface at any point are \emph{cones} \cite{Giusti}.
In case of harmonic functions or in free boundary problems, the corresponding formula implies that blow-ups are always \emph{homogeneous} \cite{CS-book,PSU}.

In the thin obstacle problem, Athanasopoulos, Caffarelli, and Salsa found that the {Almgren frequency formula}, a known monotonicity formula for harmonic functions, is still valid for solutions to the thin obstacle problem.
Such monotonicity formula states that
\[r\ \mapsto\ N(r):=\frac{r\int_{B_r(x_0)}|\nabla u|^2}{\int_{\partial B_r(x_0)}u^2}\qquad \textrm{is \, {monotone}}.\]
Thanks to this powerful tool, the blow-up sequence
\[u_r(x):=\frac{u(x_0+rx)}{\left(\ave_{\partial B_r(x_0)}u^2\right)^{1/2}}\]
converges to a \emph{homogeneous} global solution $u_0$ of degree $\mu=N(0^+)$.

Therefore, the characterization of blow-up profiles in the thin obstacle problem reduces to the characterization of homogeneous blow-up profiles.
Analyzing an eigenvalue problem on the sphere $\mathbb S^{n-1}$, and using the semi-convexity of solutions, they proved that
\[\mu<2\qquad\Longrightarrow\qquad \mu=\frac32,\]
and for $\mu=3/2$ they classified blow-ups.
Finally, using again the monotonicity property of solutions, and an appropriate boundary Harnack inequality, they established the result.

\vspace{2mm}

After the results of \cite{ACS}, the regularity of the set of regular points (i) was improved to $C^\infty$ in \cite{KPS,DS2} by using higher order boundary regularity estimates.
See Garofalo-Petrosyan~\cite{GP}, our work \cite{BFR} in collaboration with Barrios and Figalli, and the recent work of Focardi and Spadaro~\cite{FS} for a precise description and regularity results on the set of (non-regular) free boundary points satisfying (ii).

\subsection{Obstacle problems for integro-differential operators}

As explained in Section~\ref{sec3}, a natural problem that arises in probability and finance, or in physical, biological, and material sciences, is the obstacle problem for {integro-differential} operators, the simplest case being the fractional Laplacian.

In the last decade, there have been considerable efforts to extend the classical regularity theory for free boundaries of \cite{C-obst1,C-obst2} to the case of the fractional Laplacian $(-\Delta)^s$, $s\in(0,1)$.
On the one hand, this operator serves as a model case to study the regularity of the free boundary for general integro-differential operators \eqref{L}.
On the other hand, the obstacle problem for the fractional Laplacian extends at the same time the classical obstacle problem (which corresponds to the limiting case $s\to1$) and the thin obstacle problem (which corresponds to the case $s=\frac12$).

The first results in this direction were obtained by Silvestre in \cite{S-obst}, who established the almost-optimal regularity of solutions, $u\in C^{1,s-\varepsilon}$ for all $\varepsilon>0$.
The optimal $C^{1,s}$ regularity of solutions, as well as the regularity of the free boundary, were established later by Caffarelli, Salsa, and Silvestre~\cite{CSS}.

The main result of \cite{CSS} establishes that if $x_0$ is a regular point then the free boundary is $C^{1,\alpha}$ in a neighborhood of $x_0$.
More precisely, they proved that if $u$ solves the obstacle problem for the fractional Laplacian $(-\Delta)^s$ in $\R^n$, then $u\in C^{1,s}$, and for every free boundary point $x_0\in \partial\{u>\varphi\}$ we have
\begin{itemize}
\item[(i)] either \vspace{-6mm}
    \[\qquad \qquad \qquad 0<cr^{1+s}\leq \sup_{B_r(x_0)} (u-\varphi)\leq Cr^{1+s}\qquad \textrm{(regular points)}\]
\item[(ii)] or \hspace{3.3cm} $\displaystyle 0\leq \sup_{B_r(x_0)} (u-\varphi)\leq Cr^2$.
\end{itemize}
Moreover, the set of regular points (i) is an open subset of the free boundary, and it is $C^{1,\alpha}$ for some small $\alpha>0$.

Notice that the result is completely analogous to the one for the thin obstacle problem (recall that these two problems coincide if $s=\frac12$\,!).

To establish such result they found a new equivalence between the obstacle problem for the fractional Laplacian in $\R^n$ ---for \emph{every} $s\in(0,1)$---, and an appropriate \emph{thin} obstacle problem in $\R^{n+1}$.
Namely, it turns out that the fractional Laplacian $(-\Delta)^s$ can be written as a Dirichlet-to-Neumann map in $\R^{n+1}_+$ for a \emph{local} operator with a weight,
\[{\rm div}\bigl(y^{1-2s}\nabla_{x,y}\tilde u\bigr)\quad \textrm{for}\quad (x,y)\in \R^n\times\R_+;\]
see \cite{CS-ext} for more details.
When $s=1/2$, such Dirichlet-to-Neumann map is exactly the one in Remark~\ref{rem-thin}.

Thanks to such new equivalence between the obstacle problem for the fractional Laplacian and a (weighted) thin obstacle problem, they found an Almgren-type frequency formula for the obstacle problem for the fractional Laplacian in terms of such extension problem in $\R^{n+1}_+$.
Using such new monotonicity formula, they extended the regularity theory of \cite{ACS} to all $s\in(0,1)$, and also to non-zero obstacles~$\varphi$, as stated above.

\vspace{2mm}

After the results of \cite{CSS}, several new results were established concerning the structure of singular points, the higher regularity of the free boundary near regular points, or the case of operators with drift; see \cite{BFR,FS,KPS,JN,GP16}.

\vspace{2mm}

Despite all these developments in the last decade, some important problems remained open.
In particular, a very important problem that remained open after these results was the understanding of obstacle problems for more general integro-differential operators \eqref{L0}.

For the fractional Laplacian, the proofs of all known results relied very strongly on certain particular properties of such operator.
Indeed, the obstacle problem for this (nonlocal) operator is equivalent to a thin obstacle problem in $\R^{n+1}$ for a {local} operator, for which monotonicity formulas are available.

For more general nonlocal operators these tools are not available, and nothing was known about the regularity of free boundaries.
The understanding of free boundaries for more general integro-differential operators was an important problem that was completely open.

\section{Regularity theory: new results}
\label{sec5}

\subsection{Obstacle problems for general integro-differential operators}

One of my main contributions in this context is the understanding of free boundaries in obstacle problems for general integro-differential operators \cite{CRS}.
In this work we extend the results of \cite{CSS} to a much more general context, solving a long-standing open problem in the field.

Our paper \cite{CRS}, in collaboration with Caffarelli and Serra, introduces a new approach to the regularity of free boundaries in obstacle problems, and extends the results of \cite{CSS} to a general class of integro-differential operators \eqref{L}.
The main difficulty to do so was that for more general nonlocal operators $L$ there are no monotonicity formulas, while the proofs of \cite{CSS} relied strongly on such type of formulas.

Our main result in \cite{CRS} studies obstacle problems for operators \eqref{L} satisfying
\[\frac{\lambda}{|z|^{n+2s}}\leq K(z)\leq \frac{\Lambda}{|z|^{n+2s}},\qquad \textrm{with}\ K(z)\ \textrm{homogeneous},\]
and establishes that the set of regular points is open and the free boundary is $C^{1,\alpha}$ near such points.
The first assumption on $K$ is a natural uniform ellipticity assumption, and the homogeneity of $K$ is equivalent to the fact that $L$ has certain scale invariance.

More precisely, our result in \cite{CRS} establishes that, under such assumptions on $L$, if $u$ solves the obstacle problem then for every free boundary point $x_0\in \partial\{u>\varphi\}$ we have:
\begin{itemize}
\item[(i)] either \hspace{1.5cm} $\displaystyle 0<cr^{1+s}\leq \sup_{B_r(x_0)} (u-\varphi)\leq Cr^{1+s}$ \qquad \textrm{(regular points)}
\item[(ii)] or \hspace{3.7cm} $\displaystyle 0\leq \sup_{B_r(x_0)} (u-\varphi)\leq Cr^{1+s+\alpha}$,
\end{itemize}
where $\alpha>0$ is such that $1+s+\alpha<2$.
Moreover, we proved that set of {regular points} (i) is an open subset of the free boundary, and it is $C^{1,\alpha}$ for all $\alpha<s$.
Furthermore, we gave a fine description of solutions near all regular free boundary points in terms of the distance function to the free boundary.

As said before, all this was only known for the fractional Laplacian.
For more general integro-differential operators new techniques had to be developed, since one does not have any monotonicity formula.
Our proofs in \cite{CRS} are based only on very general Liouville and Harnack's type techniques, completely independent from those in~\cite{CSS}.

Let us briefly explain the global strategy of the proof.
Recall that an important difficulty is that we have no monotonicity formula, and therefore a priori blow-ups could be non-homogeneous.
As we will see, the only property we can use on blow-ups is that they are \emph{convex}.
Another difficulty is that the nonlocal operator \eqref{L} makes no sense for functions that grow too much at infinity.
Thus, we need to be very careful with the growth of functions at infinity, and the meaning of the equation as we rescale the solution and consider blow-ups.

\vspace{4mm}

\noindent \emph{Sketch of the proof}: The general argument goes as follows.
Initially, we say that a free boundary point $x_0$ is regular whenever (ii) does \emph{not} hold.
Then, we have to prove that all regular points satisfy (i), that such set is open, and that the free boundary is $C^{1,\alpha}$ near these points.

Thus, we start with a free boundary point $x_0$, and assume that (ii) does {not} hold ---otherwise there is nothing to prove.
Then, the idea is to take a blow-up sequence of the type
\[v_r(x)=\frac{(u-\varphi)(x_0+rx)}{\|u-\varphi\|_{L^\infty(B_r(x_0))}}.\]
However, we need to do it along an appropriate subsequence $r_k\to0$ so that the blow-up sequence $v_{r_k}$ (and their gradients) have a certain good growth at infinity (uniform in~$k$).
We do not want the rescaled functions $v_{r_k}$ to grow too much.
Once we do this, in the limit $r_k\to0$ we get a global solution $v_0$ to the obstacle problem, which is convex and has the following growth at infinity
\[|\nabla v_0(x)|\leq C(1+|x|^{s+\alpha}).\]
Such growth condition is very important in order to take limits $r_k\to0$ and to show that $v_0$ solves the obstacle problem.
(The fact that $v_0$ solves the obstacle problem needs to be understood in a certain generalized sense; see \cite{CRS} for more details.)

The next step is to classify global \emph{convex} solutions $v_0$ to the obstacle problem with such growth.
We need to prove that the contact set $\{v_0=0\}$ is a half-space (a priori we only know that it is convex).
For this, the first idea is to do a blow-down argument to get a new solution $\tilde v_0$, with the same growth as $v_0$, and for which the contact set is a convex \emph{cone}~$\Sigma$.
Then, we separate into two cases, depending on the size of $\Sigma$.
If $\Sigma$ has zero measure, by a Liouville theorem we show that $\tilde v_0$ would be a paraboloid, which is incompatible with the growth of $\tilde v_0$ (here we use that $s+\alpha<1$).
On the other hand, if $\Sigma$ has nonempty interior, then we prove by a dimension reduction argument (doing a blow-up at a lateral point on the cone) that $\Sigma$ must be $C^1$ outside the origin.
After that, we notice that thanks to the convexity of $\tilde v_0$ there is a cone of directional derivatives satisfying $\partial_e \tilde v_0\geq0$ in $\R^n$.
Then, using a boundary Harnack estimate in $C^1$ domains (which we prove in a separate paper \cite{RS-C1}), we show that all such derivatives have to be equal (up to multiplicative constant) in $\R^n$, and thus that $\Sigma$ must be a half-space.
Since $\Sigma$ was the blow-down of the original contact set $\{v_0=0\}$, and this set is convex, this implies that $\{v_0=0\}$ was itself a half-space.
Once we know that $\{v_0=0\}$ is a half-space, it follows that $v_0$ is a 1D solution, which can be completely classified (see \cite{RS-K}).

Once we have the classification of such blow-ups, we show that the free boundary is Lipschitz in a neighborhood of $x_0$, and $C^1$ at that point.
This is done by adapting techniques from the classical obstacle problem to the present context of nonlocal operators.
Finally, by an appropriate barrier argument we show that the regular set is open, i.e., that all points in a neighborhood of $x_0$ do \emph{not} satisfy (ii).
From this, we deduce that the free boundary is $C^1$ at every point in a neighborhood of~$x_0$, and we show that this happens with a uniform modulus of continuity around $x_0$.
Finally, using again the boundary Harnack in $C^1$ domains \cite{RS-C1}, we deduce that the free boundary is $C^{1,\alpha}$ near $x_0$.

\subsection{Application to thin obstacle problems}

Our paper \cite{CRS} not only solved an important open problem in the field, but also introduced new ideas and tools to deal with other obstacle problems when monotonicity formulas are not available.
In particular, using such methods we recently attacked the following open questions:
\begin{itemize}
\item Understanding the structure and regularity of the free boundary in thin obstacle problems for \emph{fully nonlinear} operators.
    The only known result \cite{MS} studied the regularity of solutions to this problem, but nothing was known about the free boundary.
\item Regularity of solutions and free boundaries in the thin obstacle problem with \emph{oblique} boundary condition
      \[\begin{split}
      \Delta u&=0\quad \textrm{ in}\ \R^{n+1}\cap \{x_{n+1}>0\}\\
      \min\bigl\{-\partial_{x_{n+1}} u+b\cdot\nabla u,\,u-\varphi\bigr\}&=0\quad \textrm{ on}\ \R^{n+1}\cap \{x_{n+1}=0\}.
      \end{split} \]
    Notice that such problem is equivalent to the obstacle problem
        \begin{equation}\label{drift}
        \min\bigl\{(-\Delta)^{1/2}u+b\cdot \nabla u,\,u-\varphi\bigr\}=0 \quad \textrm{in}\quad \R^n.
        \end{equation}
    The only known results in this direction studied the subcritical case $(-\Delta)^s+b\cdot \nabla$ with $s>\frac12$ \cite{GP16}, but nothing was known in the critical case $s=\frac12$.
\end{itemize}
Using the general approach introduced in \cite{CRS}, combined with important new ideas and tools needed in each of these two settings, we recently answered such open questions in \cite{RS-thin} and \cite{FR}, respectively.

The solution to these two problems presented some new and interesting features.
On the one hand, a novelty in our work \cite{RS-thin} is that we established the regularity of free boundaries without classifying blow-ups.
We succeeded in establishing a regularity result for the free boundary near regular points, but we did not prove that blow-ups are unique (a priori blow-ups could be non-unique and even non-homogeneous!).
On the other hand, in our work \cite{FR} a new phenomenon appeared: the behavior of the solution near a free boundary point $x_0$ depends on the orientation of the normal vector $\nu(x_0)$ to the free boundary (with respect to the drift $b$).
This is the first example of an obstacle problem in which this happens; see \cite{FR} for more details.

\subsection{Parabolic obstacle problems for integro-differential operators}

In collaboration with Barrios and Figalli, we studied the parabolic obstacle problem for the fractional Laplacian in \cite{BFR2}.

Despite all the developments for the elliptic problem in the last decade (described above), much less was known in the parabolic setting.
The only result was due to Caffarelli and Figalli \cite{CF}, who showed the optimal $C^{1+s}_x$ regularity of solutions in space.
However, nothing was known about the \emph{regularity of the free boundary} in the parabolic setting.

Our main theorem in \cite{BFR2} extends the results of \cite{CSS} to the parabolic setting when $s>\frac12$, and establishes the $C^{1,\alpha}$ regularity of the free boundary in $x$ and $t$ near regular points.
The result is new even in dimension $n=1$, and reads as follows.
Let us denote $Q_r(x_0,t_0)$ parabolic cylinders of size $r$ around $(x_0,t_0)$.
Then, for each free boundary point $(x_0,t_0)$, we have:
\begin{itemize}
\item[(i)] either  \hspace{1.8cm} $\displaystyle 0<c\,r^{1+s}\leq \sup_{Q_r(x_0,t_0)}(u-\varphi)\leq C\,r^{1+s}$,\vspace{2mm}
\item[(ii)] or \hspace{4cm} $\displaystyle 0\leq \sup_{Q_r(x_0,t_0)}(u-\varphi)\leq C_\varepsilon\, r^{2-\varepsilon}$ \qquad \textrm{for all}\ \,$\varepsilon>0$.
\end{itemize}
Moreover, the set of points $(x_0,t_0)$ satisfying (i) is an open subset of the free boundary and it is locally a $C^{1,\alpha}$ graph in $x$ and $t$, for some small $\alpha>0$.

Furthermore, for any point $(x_0,t_0)$ satisfying (i) there is $r>0$ such that $u\in C^{1+s}_{x,t}(Q_r(x_0,t_0))$, and we have the expansion
\begin{equation}\label{exp-parabolic}
u(x,t)-\varphi(x) = c_0\bigl((x-x_0)\cdot e+\kappa(t-t_0)\bigr)_+^{1+s}+o\bigl(|x-x_0|^{1+s+\alpha} + |t-t_0|^{1+s+\alpha}\bigr),
\end{equation}
for some $c_0>0$, $e\in \mathbb{S}^{n-1}$, and $\kappa>0$.

\vspace{2mm}

\begin{rem}[On the assumption $s>1/2$]
It is important to notice that the assumption $s>\frac12$ is necessary for this result to hold.
Indeed, the scaling of the equation in $(x,t)$ is completely different in the regimes $s>\frac12$, $s=\frac12$, and $s<\frac12$.
Since the analysis of the free boundary is always based on blow-ups, the free boundary is expected to be quite different in these three regimes.

For instance, by the examples constructed in \cite{CF}, the structure of the free boundary must be quite different when $s=\frac12$.
It was shown in \cite[Remark 3.7]{CF} that in case $s=\frac12$ there are global solutions to the parabolic obstacle problem which are homogeneous of degree $1+\beta$ for any given $\frac12\leq \beta<1$.
This means that when $s=\frac12$ there will be free boundary points satisfying neither (i) nor (ii), and there will be no ``gap'' between the homogeneities $1+s$ and 2.
This similar to what happens in the (elliptic) obstacle problem with critical drift \eqref{drift}.
\end{rem}

\vspace{2mm}

Let us briefly explain the global strategy of the proof of our result in \cite{BFR2}.
First, notice that \emph{parabolic} free boundary problems usually entail serious new difficulties with respect to their elliptic analogues.
For example, a first difficulty is that the natural scaling homogeneities are different for the equation and for the free boundary.
Namely, the equation $u_t+(-\Delta)^su=0$ is invariant under the parabolic rescalings $(x,t)\mapsto (rx,r^{2s}t)$ (recall $2s>1$), while the natural scaling of the free boundary (the one that preserves its geometry) is the hyperbolic rescaling $(x,t)\mapsto(rx,rt)$.
Thus, we are faced with a dilemma: the parabolic scaling will keep the equation but we will loose information about the free boundary, while the hyperbolic scaling will preserve the geometry of the free boundary but will loose part of the equation.
Other usual difficulties in parabolic free boundary problems are the waiting times in Harnack inequalities, and in our present setting we have the following extra difficulty: the monotonicity-type formulas of the elliptic problem do not seem to exist in the parabolic setting.
Our proof in \cite{BFR2} overcomes all these difficulties.

\vspace{4mm}

\noindent \emph{Sketch of the proof}: A first key observation is that solutions are {semiconvex} in space-time, i.e., $\partial^2_{\xi\xi}u\geq-C$ for any space-time direction $\xi\in \mathbb S^n$.
This, combined with the fact that $\partial_tu\geq0$ (a condition that comes naturally with the problem \cite{CF}), plays a very important role in our proof.
We will also use very strongly the fact that $s>1/2$.

Given any free boundary point $(x_0,t_0)$, we assume that (ii) does not hold, and we take a blow-up sequence with the parabolic rescaling.
Notice that, with such type rescaling, we expect the expansion in \eqref{exp-parabolic} to become more and more ``vertical'' in the $(x,t)$ space, so in the limit we expect the blow-up profile to be independent of~$t$ (this is what we will need to prove!).
The blow-up sequence is chosen so that the blow-up profile $u_0$ is \emph{convex} in space-time, and has \emph{subquadratic growth} at infinity.
Moreover, the origin $(0,0)$ is a free boundary point of $u_0$.

To classify blow-ups, we need to prove that any nondecreasing and convex global solution $u_0$ with subquadratic growth at infinity must be independent of $t$.
For this, we argue as follows.
First, notice that if the contact set $\{u_0=0\}$ contains a line of the form $\{(x,t)\,:\,x=\hat x\text{ for some }\hat x\in \R^n\}$, then by convexity in space-time we will have that $u_0$ is independent of $t$.
Thus, we may assume that the contact set contains no such vertical line in $(x,t)$.
But then, by monotonicity of $u_0$ in $t$, and again by convexity in $(x,t)$, we find that the contact set $\{u_0=0\}$ must be contained on a (non-vertical) half-space of the form $\{t\leq x\cdot p+c\}$, with $p\in \R^n$ and $c\geq0$; see Figure 8.

\begin{figure}[htp]
\begin{center}
\includegraphics{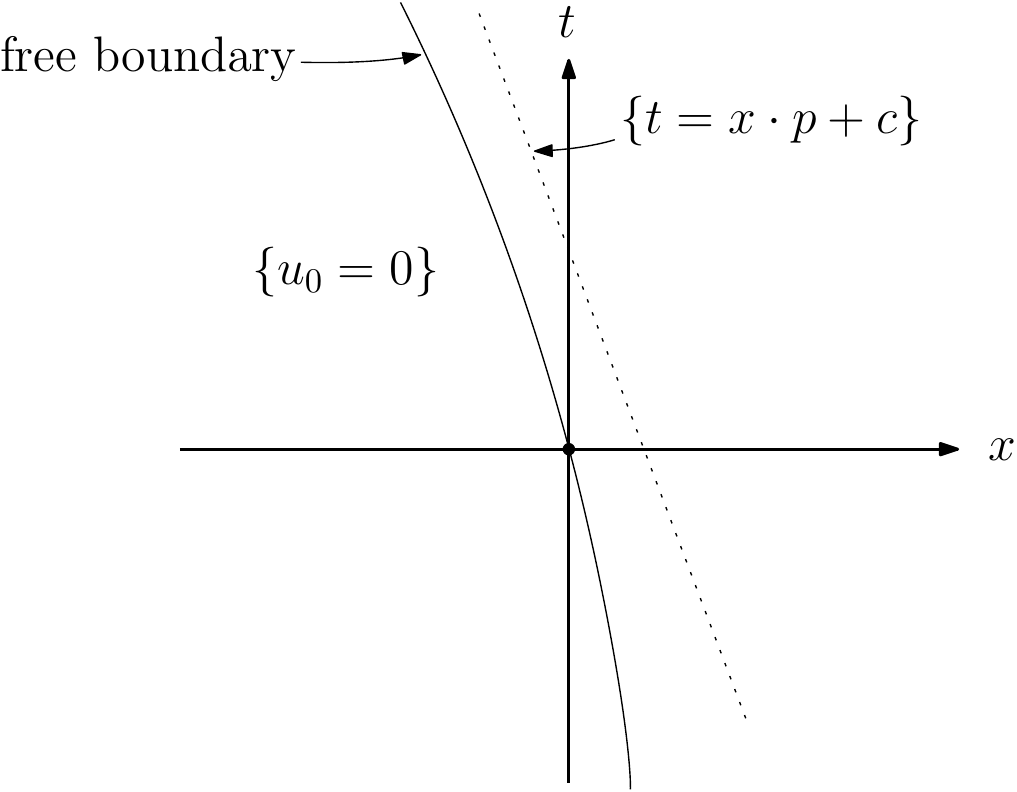}
\end{center}
\caption{The contact set $\{u_0=0\}$ must be contained in a half-space of the form $\{t\leq x\cdot p+c\}$.}
\end{figure}

We then perform a blow-down (again with the parabolic rescaling) of our global convex solution $u_0$.
Thanks to the fact that $s>\frac12$, the geometry of Figure~8 is \emph{not} preserved, and after the blow-down the set $\{t\leq x\cdot p+c\}$ becomes $\{t\leq0\}$.
Thus, we get a new global solution $U_0$, which is still nondecreasing, still convex in $(x,t)$, and still with subquadratic growth at infinity (one needs to be careful and choose the blow-down sequence to keep such growth), whose contact set $\{U_0=0\}$ is now contained in $\{t\leq0\}$.
By using the extension problem of the fractional Laplacian, combined with the inequalities satisfied by such solution $U_0$ at $t=0$ and its convexity in $x$, we show that the only possible solution is $U_0=0$ (see Lemma 3.3 in \cite{BFR2}), which in turn yields $u_0=0$.
Therefore, the only possible blow-ups $u_0(x,t)$ are those that are independent of $t$, which are completely classified in \cite{CSS} (or \cite{CRS}).

Once the blow-ups are classified, we need to transfer the information to the original solution $u$, and show that the free boundary is $C^{1,\alpha}$ in space-time near the free boundary point $(x_0,t_0)$.
For this, we first prove by standard techniques that the free boundary will be flat Lipschitz in $x$ in a neighborhood of $(x_0,t_0)$, and $C^1_x$ at that point.
With regards to regularity of the free boundary in time, we do not get that it is flat Lipschitz (this is due to the fact that we lost information when doing the blow-up with parabolic scaling), but still we prove that the free boundary is Lipschitz in $t$ (maybe with a very bad Lipschitz constant that depends on the point $(x_0,t_0)$).
Furthermore, we also show that the time derivative of $u$ is controlled by its gradient in $x$.
This is a crucial information that allows us to treat $\partial_t u$ as a lower order term, because it is controlled by $\nabla_x u$ (order 1), while the fractional Laplacian $(-\Delta)^s$ is of order $2s>1$.
Thanks to this, and using the optimal regularity in space from \cite{CF}, we can prove that $\partial_t u\in C^s_{x,t}$ and $\nabla_x u\in C^s_{x,t}$ in a neighborhood of $(x_0,t_0)$.

Then, again treating carefully $\partial_tu$ as a lower order term, and using similar ideas as in \cite{CRS}, we show that the regular set is open, and that the free boundary is $C^{1,\alpha}_x$ near $(x_0,t_0)$.
Moreover, thanks to the results of \cite{RS-C1}, we find a fine expansion of the type
\[(u-\varphi)(x,t_1)=c(x_1,t_1)d_x^{1+s}(x,t_1)+o(|x-x_1|^{1+s+\alpha})\]
for every free boundary point $(x_1,t_1)$ in a neighborhood of $(x_0,t_0)$.
Here, $d_x(\cdot,t_1)$ is the distance to the free boundary at the time slice $t=t_1$.

Finally, combining this expansion with the $C^{1,s}_t$ regularity of $u$ and the Lipschitz regularity of the free boundary in $t$, we establish the $C^{1,\alpha}$ regularity of the free boundary in space-time (maybe for a smaller $\alpha>0$).
This, combined with the expansion for $u$ at every time slice $t=t_1$, yields \eqref{exp-parabolic}, and hence the theorem.

\end{document}